\documentclass[letterpaper, 10 pt, conference]{ieeeconf}  

\IEEEoverridecommandlockouts                              
\usepackage{mathtools,amssymb,commath}
\usepackage{amsmath}
\usepackage{graphicx,import}
\usepackage{hyperref}
\hypersetup{colorlinks=true}

\usepackage{mathtools,amsthm}

\title{\LARGE \bf Momentum Control of Humanoid Robots with Series Elastic Actuators
}

\author{Gabriele Nava, Daniele Pucci and Francesco Nori$^{1}$
\thanks{*This paper was supported by the FP7 EU project CoDyCo (No. 600716 ICT 2011.2.1 Cognitive Systems and Robotics)}
\thanks{$^{1}$ All authors belong to the iCub Facility department, Istituto Italiano di Tecnologia,
        Via Morego 30, Genoa, Italy
        {\tt\small name.surname@iit.it}}%
}

\newtheorem{Remark}{\bf{Remark}}
\newtheorem{lemma}{\bf{Lemma}}

\DeclareMathOperator*{\argmin}{argmin}

\makeatletter
\newsavebox\myboxA
\newsavebox\myboxB
\newlength\mylenA
\newcommand*\xoverline[2][0.75]{%
    \sbox{\myboxA}{$\m@th#2$}%
    \setbox\myboxB\null
    \ht\myboxB=\ht\myboxA%
    \dp\myboxB=\dp\myboxA%
    \wd\myboxB=#1\wd\myboxA
    \sbox\myboxB{$\m@th\overline{\copy\myboxB}$}
    \setlength\mylenA{\the\wd\myboxA}
    \addtolength\mylenA{-\the\wd\myboxB}%
    \ifdim\wd\myboxB<\wd\myboxA%
       \rlap{\hskip 0.5\mylenA\usebox\myboxB}{\usebox\myboxA}%
    \else
        \hskip -0.5\mylenA\rlap{\usebox\myboxA}{\hskip 0.5\mylenA\usebox\myboxB}%
    \fi}
\makeatother

\begin{document}

\maketitle
\thispagestyle{empty}
\pagestyle{empty}

\begin{abstract}

Humanoid robots may require a degree of compliance at the joint level for improving efficiency, shock tolerance, and safe interaction with humans. The presence of joint elasticity, however, complexifies the design of balancing and walking controllers.
This paper proposes a  control framework for extending momentum based controllers developed for stiff actuators to the case of series elastic actuators. The key point is to consider the motor velocities as an intermediate control input, and then apply high-gain control to stabilise the desired motor velocities achieving momentum control. Simulations carried out on a model of the robot iCub verify the soundness of the proposed approach.
\end{abstract}

\section{Introduction}
\label{sec:introduction}

The paradigm  \emph{the stiffer, the better} has characterised the design of robot actuation for years. This paradigm  is well justified for position controlled industrial manipulators, where tasks usually require rapid and precise movements and very high repeatability. In the field of humanoid robotics, however, having rigid joints may be a strong limitation in terms of shock tolerance, efficiency, and safe interaction with humans \cite{Pratt1997,Tsagarakis2009}. To this purpose, compliance is usually added to the system by resorting to software techniques such as impedance or joint torque control, but the effectiveness of these control strategies might be affected by delays at any software stage \cite{Tsagarakis2009}. Another widely held solution is the redesign of robot joints by adding a spring between the load and the transmission element. The additional compliance can either vary (Variable Stiffness Actuators) or remain constant (Series Elastic Actuators) \cite{Tsagarakis2009,Parmiggiani2012,Tisi2016}. 
The presence of joint elasticity, however, complexifies the control design associated with robot manipulators. This paper presents extensions of a momentum based controller for stiff actuated  humanoid robots to the case of series elastic actuators.

In the control literature, several  techniques for controlling fixed-base robots with elastic joints have been proposed. 
Regulation tasks, for instance, can be achieved using a PD controller with on-line gravity compensation \cite{DeLuca2005}. Another possibility proposed in the last decade is to exploit the passivity properties of the overall system \cite{schaffer2007,schaffer2004}. In both cases, it is possible to prove the local (or global) stability of the closed-loop equilibrium point associated with  a set point. 
The extension of these control frameworks for addressing trajectory tracking problems may not be trivial. In the case the number of joints equals the number of actuators, a possibility is to make use of static or dynamic feedback linearization techniques \cite{spong1990control,deluca1998}. Analogously,  passivity-based controllers for regulation tasks can also be extended for addressing tracking problems \cite{schaffer2007}. These techniques, however, require the computation of joint jerks and feedforward components that can be computationally heavy for systems with a large number of degrees of freedom. 
Simplified control strategies based on \textit{singular perturbation} approaches can be designed in the case the system exhibits a two time-scale dynamics, i.e. a \textit{slow} dynamics for joint positions and a \textit{fast} dynamics for joint deformation torques~\cite{spong1990control}. The main drawback of these techniques is that they require the hypothesis of high joint stiffness.

When dealing with humanoid robots, the fixed-base assumption may be a limitation for tasks such as robot walking. An alternative solution is to use the \textit{floating base} formalism \cite{Featherstone2007}, i.e. none of the robot link has an \textit{a priori} constant position and orientation w.r.t. an inertial reference frame. However, the control problem is further complicated because the system's underactuation forbids full state feedback linearization in this case~\cite{Acosta05}. 

An effective technique for controlling floating base robots with rigid joints is the operational space control with the control task of stabilising the robot momentum. The controllers designed with this task are usually referred to as \textit{momentum-based} controllers \cite{Lee2012}. Often, the control objective for momentum based controllers  is composed of two hierarchical tasks. The primary control objective is the stabilization of robot's \textit{centroidal momentum} dynamics. This can be achieved by controlling the forces the robot exerts at contact locations \cite{Stephens2010,Herzog2014,Frontiers2015}.
The desired contact forces are obtained by  relating them to the control torques via the contact constraint equations.
Then, the secondary task exploits the redundancy of control torques (if there is any) and usually acts in the \textit{nullspace} of the primary task. The secondary task aims at the stabilization of the system \textit{zero dynamics} \cite{nava2016}. 

Extending the momentum based controllers to the case of humanoid robots powered by series elastic actuators may not be straightforward. This paper presents such an extension by assuming that the motor velocities can be assumed as an intermediate control input. This assumption allows us to retain much of the control infrastructure developed for the stiff joint case, included the gain tuning. Then, fast convergence of the  motor velocities to the desired values is achieved via feedback linearization of the motor dynamics. Simulations performed on a model of the humanoid robot iCub verify the soundness of the proposed approach.

This paper is organized as follows. Section \ref{sec:background} recalls notation, system modeling, a momentum control for robots with rigid joints, and the model for series elastic actuators. Section \ref{sec:Extension} details the proposed control strategy for controlling robots with elastic joints. Simulation results on humanoid robot iCub are presented in Section \ref{sec:results}. 
Conclusions and perspectives conclude the paper.
\section{BACKGROUND}
\label{sec:background}

\subsection{Notation}
\begin{itemize}
\item $\mathcal{I}$ denotes an inertial frame, with its $z$ axis pointing against the gravity. The constant $g$ denotes the norm of the gravitational acceleration.
\item Given a matrix $A \in \mathbb{R}^{m \times n}$, we denote with $A^{\dagger}\in \mathbb{R}^{n \times m}$ its Moore Penrose pseudoinverse. 
\item $e_i \in  \mathbb{R}^m$ is the canonical vector, consisting of all zeros but the $i$-th component that is equal to one.
\item We denote with $m$ the total mass of the robot.
\end{itemize}

\subsection{Modelling with rigid transmissions}
The robot is modelled as a multi-body system composed of $n + 1$ rigid bodies, called links, connected by $n$ joints with one degree of freedom each. We also assume that
none of the links has an \emph{a priori} constant pose with respect to an inertial frame, i.e. the system is \emph{free floating}.

The robot configuration space is the Lie group $\mathbb{Q} = \mathbb{R}^3 \times SO{(3)} \times \mathbb{R}^n$  and it is characterized by the \emph{pose} 
(position and orientation) of a \emph{base frame} attached to a robot's link, and the joint positions.
An element $q \in \mathbb{Q}$ can be defined as the following triplet: $q = (\prescript{\mathcal{I}}{}p_{\mathcal{B}}, \prescript{\mathcal{I}}{}R_{\mathcal{B}}, s)$ 
where $\prescript{\mathcal{I}}{}p_{\mathcal{B}} \in \mathbb{R}^3$ denotes the position of the base frame with respect to the inertial frame,
$\prescript{\mathcal{I}}{}R_{\mathcal{B}} \in \mathbb{R}^{3\times3}$ is a rotation matrix representing the orientation of the base frame, and $s \in \mathbb{R}^n$ is the 
joint configuration characterising the \emph{shape} of the robot. 
The velocity of the multi-body system can be characterized by the set $\mathbb{V} = \mathbb{R}^3 \times \mathbb{R}^3 \times \mathbb{R}^n$.
An element of $\mathbb{V}$ is a triplet $\nu = ( ^\mathcal{I}\dot{ p}_{\mathcal{B}},^\mathcal{I}\omega_{\mathcal{B}},\dot{s}) = (\text{v}_{\mathcal{B}}, \dot{s})$,
where $^\mathcal{I}\omega_{\mathcal{B}}$ is the angular velocity of the base frame expressed w.r.t. the inertial frame, 
i.e. $^\mathcal{I}\dot{R}_{\mathcal{B}} = S(^\mathcal{I}\omega_{\mathcal{B}})^\mathcal{I}{R}_{\mathcal{B}}$.
A more detailed description of the floating base model is provided in \cite{traversaro2017}.

We assume that the robot interacts with the environment by exchanging $n_c$ distinct wrenches. The equations of motion of the multi-body system can be described
applying the Euler-Poincar\'e formalism \cite[Ch. 13.5]{Marsden2010}:
\begin{align}
\label{eq:system}
{M}(q)\dot{{\nu}} + {C}(q, {\nu}){\nu} + {G}(q) =  B \tau + \sum_{k = 1}^{n_c} {J}^\top_{\mathcal{C}_k} f_k
\end{align}
where ${M} \in \mathbb{R}^{n+6 \times n+6}$ is the mass matrix, ${C} \in \mathbb{R}^{n+6 \times n+6}$ is the Coriolis matrix, ${G} \in \mathbb{R}^{n+6}$ is the gravity 
term, $B = (0_{n\times 6} , 1_n)^\top$ is a selector matrix, $\tau \in \mathbb{R}^{n}$ is a vector representing the internal actuation torques, and 
$f_k \in \mathbb{R}^{6}$ denotes an external wrench applied by the environment to the link of the $k$-th contact. The Jacobian 
${J}_{\mathcal{C}_k} = {J}_{\mathcal{C}_k}(q)$ is the map between the robot's velocity ${\nu}$ and the linear and angular velocity at the $k$-th contact link.

As described in \cite[Sec. 5]{traversaro2017}, it is possible to apply a coordinate transformation in the state space $(q,{\nu})$ that transforms the system dynamics~\eqref{eq:system} into a new form where
the mass matrix is block diagonal, thus decoupling joint and base frame accelerations. Also, in this new set of coordinates,  the first six rows of Eq. \eqref{eq:system} are the \emph{centroidal dynamics}\footnote{In the specialized literature, the term \emph{centroidal dynamics} 
is used to indicate the rate of change of the robot's momentum expressed at the center-of-mass, which then equals the summation of all external wrenches acting on the 
multi-body system \cite{orin2013}.}. As an abuse of notation, we assume that system \eqref{eq:system} has been transformed into this new set of coordinates, i.e. 
\begin{IEEEeqnarray}{RCL}
\label{centrTrans}
M(q) &=& \begin{bmatrix} {M}_b(q) & 0_{6\times n} \\ 0_{n\times 6} & {M}_j(q) \end{bmatrix}, \quad
H    = M_b \text{v}_{\mathcal{B}},
\end{IEEEeqnarray}
where ${M}_b \in \mathbb{R}^{6\times 6}, {M}_j \in \mathbb{R}^{n\times n}$, ${H}:=(H_L,H_\omega)\in \mathbb{R}^6$ is the robot centroidal momentum, and $H_L, H_\omega \in \mathbb{R}^3$ are the linear and angular momentum at the center of mass, respectively.
                           
Lastly, it is assumed that a set of holonomic constraints acts on System \eqref{eq:system}. These holonomic constraints are of the form $c(q) = 0$, and may represent,
for instance, a frame having a constant pose w.r.t. the inertial frame.
In the case where this frame corresponds to the location at which a rigid contact occurs on a link, we represent the holonomic constraint as
$
{J}_{\mathcal{C}_k}(q) {\nu} = 0.$
Hence, the holonomic constraints associated with all the rigid contacts can be represented as
\begin{IEEEeqnarray}{RCL}
\label{eqn:constraintsAll}
{J}(q) {\nu} {=} 
\begin{bmatrix}{J}_{\mathcal{C}_1}(q) \\ \cdots \\ {J}_{\mathcal{C}_{n_c}}(q)  \end{bmatrix}{\nu}  {=} 
\begin{bmatrix} J_b & J_j  \end{bmatrix}{\nu}  
&=& J_b {\text{v}}_{\mathcal{B}}+ J_j \dot{s} = 0,
\IEEEeqnarraynumspace
\end{IEEEeqnarray}
with $J_b \in \mathbb{R}^{6n_c \times 6},J_j \in \mathbb{R}^{6n_c \times n} $.
The base frame velocity is denoted by $\text{v}_{\mathcal{B}} \in \mathbb{R}^6$, which in the new coordinates yielding a block-diagonal mass matrix is given by $\text{v}_{\mathcal{B}} = (\dot{p}_c,\omega_o)$, where $\dot{p}_c \in \mathbb{R}^3$ is the velocity of the system's center of mass ${p}_c \in \mathbb{R}^3$, and $\omega_o \in \mathbb{R}^3$ is the so-called system's \emph{average angular velocity}.
By differentiating the kinematic constraint Eq. \eqref{eqn:constraintsAll}, one obtains
\begin{equation}
    \label{eq:constraints_acc}
   J\dot{\nu}+\dot{J}\nu = J_b \dot{\text{v}}_{\mathcal{B}}+ J_j \ddot{s} + \dot{J}_b {\text{v}}_{\mathcal{B}}+ \dot{J}_j \dot{s} = 0.
\end{equation}

In view of \eqref{eq:system}--\eqref{centrTrans}, the equations of motion along the constraints \eqref{eq:constraints_acc} are given by:
\begin{IEEEeqnarray}{LCL}
    \label{eq:systemCentroidal}
    M_b\dot{\text{v}}_\mathcal{B}   & = &  J_b^{\top}f - h_b                                              \IEEEyessubnumber \label{eq:floatingCentroidal} \\
    M_j\ddot{s}    & = &  J_j^{\top}f - h_j  + \tau                                                       \IEEEyessubnumber \label{eq:jointsCentroidal}   
\end{IEEEeqnarray} 
where we define $h:={C}(q, {\nu}){\nu} + {G}(q)  \in \mathbb{R}^{n+6}$ and its partition $h = (h_b,h_j), h_b\in\mathbb{R}^6,  h_j\in\mathbb{R}^n$. $f:=(f_1,\cdots,f_{n_c}) \in \mathbb{R}^{6n_c} $ are the set of contact forces -- i.e. Lagrange multipliers -- making Eq.~\eqref{eq:constraints_acc} satisfied.

\subsection{Balancing control with robot rigid transmissions}
\label{subsec:rigidControl}

We recall here  the \emph{momentum-based} control strategy implemented on the iCub humanoid robot~\cite{nava2016,pucciTuning2016,pucciBalancing2016}. 
The control objective is the stabilization of a desired robot momentum and the stability of the associated zero dynamics. 

\subsubsection{Momentum control}

Recall that the rate-of-change of the robot momentum equals the net external wrench acting on the robot, which in the present case reduces to the contact
wrench $f$ plus the gravity wrenches. In view of Eq. \eqref{centrTrans}, the rate-of-change of the robot momentum can be expressed as:
\begin{IEEEeqnarray}{RCCCL}
\label{hDot}
 \frac{\dif }{\dif t}(M_b {\text{v}_\mathcal{B}}) &=& \dot{H}(f) &=& J_b^{\top}f - mge_3,
\end{IEEEeqnarray}
where $e_3 \in \mathbb{R}^6$.
Let $H^d \in \mathbb{R}^6 $ denote the desired robot momentum, and  $\tilde{H} = H - H^d$  the momentum error. Assuming  that the contact wrenches $f$ can be chosen at will, then we choose $f$ such that \cite{nava2016}:
\begin{IEEEeqnarray}{RCL}
    \label{hDotDes}
    \dot{H}(f) &=& 
    \dot{H}^* := \dot{H}^d - K_p \tilde{H} - K_i I_{\tilde{H}}    
    \IEEEeqnarraynumspace  \IEEEyessubnumber \label{eq:dotH_f}  \\
    \dot{I}_{\tilde{H}} &=&
                     \begin{bmatrix} 
                     {J}_{G}^L(s) \\ 
                     {J}_{G}^{\omega}(s^d)
                     \end{bmatrix}\dot{s} \IEEEyessubnumber  \label{eq:IhTilde} 
\end{IEEEeqnarray}
$K_p, K_i {\in} \mathbb{R}^{6\times 6}$ two symmetric, positive definite matrices and 
\begin{IEEEeqnarray}{RCL}
\label{eqn:jacobian} 
\bar{J}_G(s) &{:=}& - M_bJ^{\dagger}_bJ_j =  \begin{bmatrix} {J}_{G}^L(s) \\ {J}_{G}^{\omega}(s)\end{bmatrix} \in \mathbb{R}^{6\times n}, {J}_{G}^L, {J}_{G}^{\omega}\in \mathbb{R}^{3\times n}.\IEEEeqnarraynumspace \nonumber
\end{IEEEeqnarray} 
If $n_c > 1$, there are infinite contact wrenches $f$ that satisfy Eq.~\eqref{eq:dotH_f}.
We parametrize the set of solutions $f$ to \eqref{eq:dotH_f} as:
\begin{equation}
    \label{eq:forces}
    f = f_1 + N_{b}f_0
\end{equation}
with $f_1 =  J_b^{\top\dagger} \left(\dot{H}^*+ mg e_3\right)$, $N_b \in \mathbb{R}^{6n_c \times 6n_c}$ the projector into the null space of $J_b^{\top}$, and $f_0\in \mathbb{R}^{6n_c}$ the wrench redundancy that does not influence $\dot{H}(f) =
    \dot{H}^*$.
To determine the control torques that instantaneously realize the contact wrenches given by \eqref{eq:forces}, we use the dynamic 
equations \eqref{eq:system} along with the constraints \eqref{eq:constraints_acc}, which yields:
\begin{equation}
    \label{eq:torques}
    \tau = \Lambda^\dagger (JM^{-1}(h - J^\top f) - \dot{J}\nu) + N_\Lambda \tau_0
\end{equation}
with $\Lambda = {J_j}{M_j}^{-1} \in \mathbb{R}^{6n_c\times n}$,  $N_\Lambda \in \mathbb{R}^{n\times n}$  the projector onto the nullspace
of $\Lambda$,  and $\tau_0 \in \mathbb{R}^n$  a free variable.

\subsubsection{Stability of the Zero Dynamics}
The stability of the zero dynamics is attempted by means of a so called \emph{postural task}, which exploits the free variable $\tau_0$ in \eqref{eq:torques}.
A choice of the postural task that ensures the stability of the zero dynamics on one foot is \cite{nava2016}:
\begin{IEEEeqnarray}{lCr}
\label{posturalNew}
    \tau_0 &=& h_j - J_j^\top f + u_0
\end{IEEEeqnarray}
where 
$u_0 := -K^j_{p}N_\Lambda M_j(s-s^d) -K^j_{d}N_\Lambda M_j\dot{s}$, and
$K^{j}_p \in \mathbb{R}^{n \times n}$ and $K^{j}_d \in \mathbb{R}^{n \times n}$ two symmetric, positive definite matrices. An interesting property of the closed loop system~\eqref{eq:system}--\eqref{eq:forces}--\eqref{eq:torques}--\eqref{posturalNew} is recalled in the following Lemma.

\begin{lemma}[~\cite{pucciTuning2016}]
\label{lemmaf0}
Assume that $\Lambda$ is full row rank. Then, the closed loop joint space dynamics $\ddot{s}$ does not depend upon the wrench redundancy $f_0$.
\end{lemma}

This result is a consequence of the postural control choice~\eqref{posturalNew} and it is of some interest: it means that the closed loop joint dynamics depends on the total rate-of-change of the momentum, i.e. $\dot{H}$, but not on the different forces generating it. Hence, any choice of the redundancy $f_0$ does not influence the joint dynamics $\ddot{s}$, and we  exploit it to minimize the joint torques $\tau$ in Eq.~\eqref{eq:torques}.

In the language of \emph{Optimization Theory}, we can rewrite the  control strategy as a single optimisation problem as follows:
\begin{IEEEeqnarray}{RCL}
	\IEEEyesnumber
	\label{inputTorquesSOT}
	f^* &=& \argmin_{f}  |\tau^*(f)|^2 \IEEEyessubnumber \label{inputTorquesSOTMinTau}  \\
		   &s.t.& \nonumber \\
		   &&Cf < b \IEEEyessubnumber  \label{frictionCones} \\
		   &&\dot{H}(f) = \dot{H}^*  \IEEEyessubnumber \\
		   &&\tau^*(f) = \argmin_{\tau}  |\tau - \tau_0(f)|^2 \IEEEyessubnumber	\label{optPost} 
  \\
		   	&& \quad s.t.  \nonumber \\
		   	&& \quad \quad \ \dot{J}(q,\nu)\nu + J(q)\dot{\nu} = 0
		    \IEEEyessubnumber 	\label{constraintsRigid} \\
		   	&& \quad \quad \ \dot{\nu} = M^{-1}(B\tau+J^\top f {-} h) \IEEEyessubnumber \\
		   && \quad \quad \ 	\tau_0 = 
		   h_j - J_j^\top f + u_0.		    \IEEEyessubnumber
\end{IEEEeqnarray}
The constraints~\eqref{frictionCones} ensure the satisfaction of friction cones, normal contact surface forces, and center-of-pressure constraints. The control torques are then given by $\tau {=} \tau^*(f^*)$.

\subsection{Series Elastic Actuator actuation case}

The previous section has presented a balancing control assuming that the joint torque $\tau$ can be assumed as control input. This is the case, for instance, of a torque controlled robot where the motors are rigidly connected to the joints, eventually by means of harmonic drives. In this section, we present the extension of the model~\eqref{eq:system} when the interfaces between the motors and the joints are elastic elements, namely, the robot is powered by series elastic actuators. To this purpose, we make the following assumptions.

\begin{figure}[t]
   \centering
   \includegraphics[width=\columnwidth]{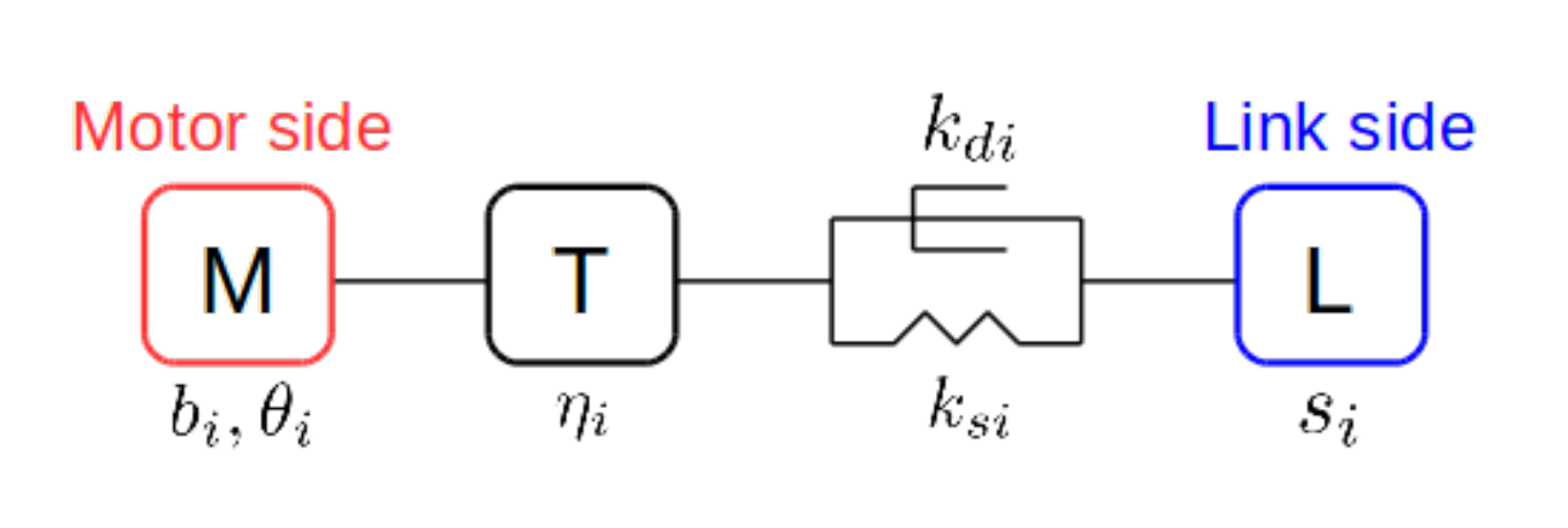}
   \caption{Block diagram of the series elastic actuator.}
   \label{MTU}
\end{figure}

\begin{itemize}
    \item The angular motor kinetic energy is due to  its own spinning only, and the center of masses of each motor is along the motor axis of rotation;
    \item Both  stiffness and damping of the series elastic actuators can be considered linear versus its relative displacement absolute value and rate-of-change.
    \item All motors are rigidly connected to the transmission element.
\end{itemize}
Figure~\ref{MTU} depicts a simple block diagram of the series elastic actuator assumed to power the robot. In this picture, $\theta_i$ is the $i$-th motor position,  $s_i$  the $i$-th link position, $\eta_i$ is the transmission ratio, $b_i$ is the motor inertia, and $k_{si}$ and $k_{di}$ are the $i$-th link (torsional) stiffness and damping. 

Under the above assumptions, the model~\eqref{eq:systemCentroidal} can be extended by adding the dynamics of the motor angles $\theta = \begin{bmatrix} \theta_1 &  ... & \theta_n \end{bmatrix}^{\top} \in \mathbb{R}^n$. Then (see~ \cite{schaffer2004},\cite{deluca1998},\cite{spong1990control} for details):
\begin{IEEEeqnarray}{LCL}
    \label{eq:elasticSystem}
    M_b\dot{\text{v}}_\mathcal{B}   & = &  J_b^{\top}f - h_b                                              \IEEEyessubnumber \label{eq:floating}  \\
    M_j\ddot{s}    & = &  J_j^{\top}f - h_j  + \tau       \IEEEyessubnumber \label{eq:joints}    \\
    I_m\ddot{\theta} & = & \tau_m - \Gamma\tau             \IEEEyessubnumber \label{eq:motors} 
\end{IEEEeqnarray}
where 
\begin{IEEEeqnarray}{LCL}
    \label{eq:alpha}
    \tau := K_S(\Gamma\theta-s) + K_D(\Gamma\dot{\theta}-\dot{s}) \in \mathbb{R}^n
\end{IEEEeqnarray}
represents the coupling  between the joints dynamics and the motors dynamics. The positive definite diagonal matrices describing joints stiffness and damping are $K_S = \text{diag}(k_{si})$ and $K_D = \text{diag}(k_{di})$, respectively. The positive definite matrix $\Gamma = \text{diag}(\eta_i)\in \mathbb{R}^{n \times n}$ is a \emph{transmission matrix}, and $I_m = \text{diag}(b_i) \in \mathbb{R}^{n \times n}$ is the motor inertia matrix. The control input is given by the motor torques $\tau_m \in \mathbb{R}^n$. The new system configuration space is the Lie Group \[\overline{\mathbb{Q}} = \mathbb{R}^3 \times SO(3) \times \mathbb{R}^{2n},\] and a system configuration is represented by the quadruple $\overline{q} = (\prescript{\mathcal{I}}{}p_{\mathcal{B}},\prescript{\mathcal{I}}{}R_{\mathcal{B}},s,\theta)$. The velocity is then represented by the set $\overline{\mathbb{V}} = \mathbb{R}^3 \times \mathbb{R}^3 \times \mathbb{R}^{2n}$, and an element of $\overline{\mathbb{V}}$ is then $\overline{\nu} = (^\mathcal{I}\dot{p}_{\mathcal{B}},^\mathcal{I}\omega_{\mathcal{B}},\dot{s},\dot{\theta})$. The constraint equation~\eqref{eq:constraints_acc} remains invariant, namely, it is not affected by the addition of the motor dynamics. In fact, this equation represents the fact that the feet acceleration, expressed in terms of $q$, is equal to zero for all the time.

\begin{Remark}
\label{remark1}
We assume that the matrix $M_j$ in Eq.~\eqref{centrTrans} does not contain terms due to the so-called \emph{motor reflected inertia}. More precisely, the equations of motion for the rigid actuation case can be deduced by imposing $K_S \xrightarrow[]{} \infty$ in Eq. \eqref{eq:elasticSystem}, which implies $s \xrightarrow[]{} \Gamma\theta$, $\dot{s} \xrightarrow[]{} \Gamma \dot{\theta}$, $\ddot{s} \xrightarrow[]{} \Gamma \ddot{\theta}$. Then,  by summing up \eqref{eq:joints}-\eqref{eq:motors}, and by multiplying \eqref{eq:motors} times $\Gamma^{-1}$, one has:
\emph{\begin{IEEEeqnarray}{LCL}
    \label{eq:equivRigidSystem}
    M_b\dot{\text{v}}_\mathcal{B}   & = &  J_b^{\top}f - h_b                                                 \IEEEyessubnumber \label{eq:equivFloating}  \\
    (M_j+\Gamma^{-1}I_m\Gamma^{-1})\ddot{s}    & = &  J_j^{\top}f - h_j  + \Gamma^{-1}\tau_m.  \IEEEyessubnumber \label{eq:equivJoints}    
\end{IEEEeqnarray}}
These equations of motion characterise the system evolution in the case of rigid transmissions, and the term $\Gamma^{-1}I_m\Gamma^{-1}$ is called \emph{motor reflected inertia}. Hence, we assume that $M_j$ does not take this term into account. 
\end{Remark}
\section{Control Design}
\label{sec:Extension}
The control of system~\eqref{eq:elasticSystem} for robot balancing purposes may not be straightforward. In fact, assuming that the control objective is still the asymptotic stabilisation of the momentum error $\tilde{H}$, then its rate-of-change is no longer influenced by the system control input, namely, the motor torque $\tau_m$. More precisely, the momentum equation~\eqref{hDot} still holds, i.e.
\begin{IEEEeqnarray}{RCL}
  \dot{H}(f) &=& J_b^{\top}f - mge_3. \nonumber
\end{IEEEeqnarray}
By using the (feet zero-acceleration) constraint~\eqref{eq:constraints_acc}, i.e.
\begin{equation}
   J_b \dot{\text{v}}_{\mathcal{B}}+ J_j \ddot{s} + \dot{J}_b {\text{v}}_{\mathcal{B}}+ \dot{J}_j \dot{s} = 0, \nonumber
\end{equation}
 with the robot acceleration deduced from~\eqref{eq:floating}-\eqref{eq:joints}, i.e.
\begin{IEEEeqnarray}{LCL}
   \dot{\text{v}}_{\mathcal{B}}   & = &   M^{-1}_b \left( J_b^{\top}f - h_b \right)
    \nonumber \\
   \ddot{s}    & = &  M^{-1}_j \left( J_j^{\top}f - h_j  + \tau(\theta,\dot{\theta},s,\dot{s})  \right)   
    \nonumber
\end{IEEEeqnarray}
one observes that we can no longer  relate the contact force~$f$ to the system input $\tau_m$. Consequently, the momentum rate-of-change $\dot{H}$ is no longer instantaneously influenced by the system input $\tau_m$.

For the same reason, the control algorithm recalled in section~\ref{subsec:rigidControl} is no longer applicable in the case of robots powered by series elastic actuators. In fact, even if one chooses a contact force $f$ to achieve $\dot{H}(f) = \dot{H}^*$, then this force can no longer be achieved by~\eqref{eq:torques}, i.e.
\begin{IEEEeqnarray}{RCL}
    \tau = \Lambda^\dagger (JM^{-1}(h - J^\top f) - \dot{J}\nu) + N_\Lambda \tau_0, \nonumber
\end{IEEEeqnarray}
since $\tau$ does not depend upon the control input $\tau_m$, but only on the system state --see Eq.~\eqref{eq:alpha}.

In the language of \emph{Control Theory}, we would say that the output $H$ has a \emph{relative degree} equal to two: the second order time derivative of $H$ can be instantaneously influenced by the control input $\tau_m$. Then, one may attempt at the control of $H$ by performing feedback linearisation of $\ddot{H}$, and then choose the input redundancy for postural control. 

Analogously, one may consider $\dot{f}$ as fictitious control input in the dynamics of $\ddot{H}$, i.e.
\begin{IEEEeqnarray}{RCL}
  \ddot{H}(\dot{f}) &=& \dot{J}_b^{\top}f + J_b^{\top}\dot{f}, \nonumber
\end{IEEEeqnarray}
and then exploit the time derivative of~\eqref{eq:torques} with \eqref{eq:motors} to impose a rate-of-change of the contact force $f$.

The main drawbacks of these strategies -- which are based on pure-feedback linearisation of an output with relative degree higher than one -- are the following:
\begin{itemize}
    \item they need of feedforward terms seldom precisely known in practice, such as $\dot{M}$, $\ddot{J}$, etc;
    \item when series-elastic-actuators are introduced to substitute some, or all, rigid transmission mechanisms, they usually need a new time-consuming and specific  gain tuning procedure;
    \item they do not leverage the (usually) high frequency and reliable low-level motor velocity control.
\end{itemize}
We then follow a different route to deal with series elastic actuators that aims -- when the  motors are equipped with it -- at exploiting the low level motor velocity control. In addition, the proposed strategy is shown to be robust against the feedforward terms usually needed by feedback linearisation, and can also exploit the gain tuning procedure developed for the rigid actuation case.

\subsection{The balancing control  for series elastic actuators}
Eq.~\eqref{eq:motors}  points out that the motor dynamics $\ddot{\theta}$ is fully actuated. Then, any desired motor velocity $\dot{\theta}_d$ can be stabilised with any desired, \emph{small}, settling time. This in turn implies that the motor velocity can be assumed as a virtual control input in the dynamics~\eqref{eq:joints}. 
In the language of Automatic Control, 
assuming $\dot{\theta}$ as control variable  is a typical \emph{backstepping} assumption. Then, the production of  the motor torques associated with the desired motor velocities 
can be achieved via classical nonlinear control techniques~\cite[p. 589]{kh02} or high-gain control. We later detail an implementation of the latter for obtaining the aforementioned motor torques and perform simulation results.

More precisely, we consider $\beta=K_D\Gamma\dot{\theta}$ as a fictitious control input in the joints dynamics~\eqref{eq:joints}. Then, one has:
\begin{IEEEeqnarray}{LCL}
    \label{newJointDynamics}
    M_s\ddot{s}    & = &  J_j^{\top}f - \bar{h}_j  + \beta,
\end{IEEEeqnarray}
with 
\begin{IEEEeqnarray}{LCL}
    \label{newJointDynamicsElements}
     \bar{h}_j &=& = {h}_j - p   \IEEEyessubnumber \label{hjBar} \\
     p &=& K_S(\Gamma\theta-s) - K_D\dot{s} \IEEEyessubnumber \label{p}.
\end{IEEEeqnarray}
Observe that~\eqref{newJointDynamics} concides with~\eqref{eq:jointsCentroidal} by substituting $h_j$ with $\bar{h}_j$ and $\tau$ with $\beta$. In light of this, the control input $\beta$ achieving balancing control -- with the same objectives detailed in section~\ref{subsec:rigidControl} -- is achieved by solving the optimisation problem~\eqref{inputTorquesSOT} with the aforemetioned substitutions, i.e.
\begin{IEEEeqnarray}{RCL}
	\IEEEyesnumber
	\label{inputTorquesSOTu}
	f^* &=& \argmin_{f}  |\beta^*(f)|^2 \IEEEyessubnumber \label{inputTorquesMinVel}  \\
		   &s.t.& \nonumber \\
		   &&Cf < b \IEEEyessubnumber   \\
		   &&\dot{H}(f) = \dot{H}^*  \IEEEyessubnumber \\
		   &&\beta^*(f) = \argmin_{\beta}  |\beta - \beta_0(f)|^2 \IEEEyessubnumber
  \\
		   	&& \quad s.t.  \nonumber \\
		   	&& \quad \quad \ \dot{J}(q,\nu)\nu + J(q)\dot{\nu} = 0
		    \IEEEyessubnumber 	 \\
		   	&& \quad \quad \ \dot{\nu} = M^{-1}(B\beta+J^\top f {-} h) \IEEEyessubnumber \\
		   && \quad \quad \ 	\beta_0 = 
		   \bar{h}_j - J_j^\top f + u_0.		    \IEEEyessubnumber
\end{IEEEeqnarray}
with $h := (h_b,\bar{h}_j)$, and \[u_0 := -K^j_{p}N_\Lambda M_j(s-s^d) -K^j_{d}N_\Lambda M_j\dot{s}.\] 

The optimisation problem~\eqref{inputTorquesSOTu} points out that the redundancy of the contact forces in achieving $\dot{H}(f) = \dot{H}^*$ is no longer exploited to minimise the joint torques, but rather the torques induced by the motor velocities -- compare  Eqs.~\eqref{inputTorquesSOTMinTau} and~\eqref{inputTorquesMinVel}. This basically means that the solution to the above  problem tends to minimise the desired motor velocities.

Once the optimisation problem~\eqref{inputTorquesSOTu} is solved, at each time instant the (desired) motor velocities are given by 
\begin{IEEEeqnarray}{RCL}
	\IEEEyesnumber
	\label{inputDesMotVel}
	\dot{\theta}_d = \Gamma^{-1} K^{-1}_D \beta^*(f^*).
\end{IEEEeqnarray}
Now, if the series-elastic actuators provide the user with a velocity control interface,  one can send as desired values to this interface the velocities~\eqref{inputDesMotVel}. On the other hand, if the series-elastic actuators provide as interface the motor torques~$\tau_m$ (often related to the motor $PWM$), then one can apply high-gain control on the dynamics~\eqref{eq:motors} to stabilise the motor velocity~\eqref{inputDesMotVel}. In particular, motor velocity control via the motor torque $\tau_m$ can be achieved by:
\begin{IEEEeqnarray}{RCL}
    \label{eq:tauMotors}
    \tau_m & = & -I_mK_{m}(\dot{\theta}-\dot\theta_d) + \Gamma\tau
\end{IEEEeqnarray}
with $K_{m} \in \mathbb{R}^{n\times n}$ a positive diagonal matrix. The closed-loop motors dynamics is then given by:
\begin{IEEEeqnarray}{RCL}
    \label{eq:controlledMotors}
    \ddot{\theta} & = & -K_{m}(\dot{\theta}-\dot\theta_d),
\end{IEEEeqnarray}
which implies that the motor velocity tracking error stays relatively \emph{small} for relatively  high gains $K_{m}$. 
It is important to observe that the control law~\eqref{eq:tauMotors} misses the feed-forward component  $\ddot\theta_d$, and it is not deduced by the application of pure backstepping techniques~\cite[p. 589]{kh02}. All these missing elements represent a robustness test for the controller presented in this paper.


\subsection{The mixed actuation case: stiff and elastic actuators}
The framework presented above may be useful when the humanoid robot is powered by  stiff and elastic actuators. In this case, one can still solve the optimisation problem~\eqref{inputTorquesSOTu} by properly defining the vectors $\bar{h}_j$ and $p$ in Eq.~\eqref{newJointDynamicsElements}. 

To provide the reader with an example of such a mixed actuation case, assume that the robot possesses $m_1$ stiff actuators -- for which the associated joint torques $\tau_{m_1} \in \mathbb{R}^{m_1}$ can be considered as control inputs -- and $m_2$ series-elastic actuators, with $n = m_1 + m_2$. Assume that in the serialisation of the joint angles $s = (s_{m_1},s_{m_2})$, the first $m_1$ joints are powered by  stiff actuators, and the remaining $m_2$ by  elastic ones. Then, the optimisation problem~\eqref{inputTorquesSOTu} can be solved with 
\begin{IEEEeqnarray}{LCL}
\beta &:=& 
\begin{pmatrix}
    \tau_{m_1} \\
    K_D \Gamma \dot{\theta}
\end{pmatrix} \\
     p &=& 
\begin{pmatrix}
    0\\
    K_S(\Gamma\theta-s_{m_2}) - K_D\dot{s}_{m_2}
\end{pmatrix}
     \IEEEyessubnumber 
\end{IEEEeqnarray}
where now $\theta \in \mathbb{R}^{m_2}$. 
\section{Simulation results}
\label{sec:results}

We test the proposed control solution by using a model of the humanoid robot iCub~\cite{Metta20101125} with 25 degrees-of-freedom (DoFs). The inertia values $b_i$ are obtained from the motor datasheets, and their order of magnitude is around $10^{-5}$ $[kg \ m^2]$. Realistic stiffness and damping values for the series elastic actuators are obtained from previous work on the design of SEA for iCub legs \cite{Parmiggiani2012,Tisi2016}. The stiffness value used for the experiments is $350$ $[\frac{Nm}{rad}]$ for all the joints while damping coefficient is $0.25$ $[\frac{Nms}{rad}]$. The transmission ratio is the same for all joints and set equal to $\eta_i = \frac{1}{100}$.

\subsection{Simulation Environment}

The control algorithm is tested by performing simulations in  the  MATLAB environment. In particular, recall that $\overline{q} \in \overline{\mathbb{Q}} = \mathbb{R}^3 \times SO(3) \times \mathbb{R}^{2n}$: it is then necessary to choose a representation for the Special Orthogonal Group $SO(3)$. We use quaternion parametrization instead of Euler angles because it does not introduce artificial singularities. The resulting state space and its time derivative are given by:
\begin{IEEEeqnarray}{RCL}
    \chi & := & \begin{bmatrix} ^\mathcal{I}{}p_{\mathcal{B}} & \mathcal{Q} & s & \theta & ^\mathcal{I}{}\dot{p}_{\mathcal{B}} & ^\mathcal{I}{}\omega_{\mathcal{B}} & \dot{s} & \dot{\theta} \end{bmatrix}^{\top} \\
    \dot{\chi} & = & \begin{bmatrix} ^\mathcal{I}{}\dot{p_{\mathcal{B}}} & \dot{\mathcal{Q}} & \dot{s} & \dot{\theta} & \dot{\overline{\nu}} \end{bmatrix}^{\top}. \label{eq:dotchi}
\end{IEEEeqnarray}
The system evolution is obtained by integrating Eqs.~\eqref{eq:dotchi}-\eqref{eq:elasticSystem} using MATLAB numerical integrator \textit{ode15s}.
The numerical integration, as well as the constraint $|\mathcal{Q}|=1$, are enforced during the integration phase using additional integration terms \cite{Gros2015}.

\subsection{Tracking Performances}
\label{subsec:tracking}

\begin{figure}[t!]
   \centering
   \includegraphics[width=\columnwidth]{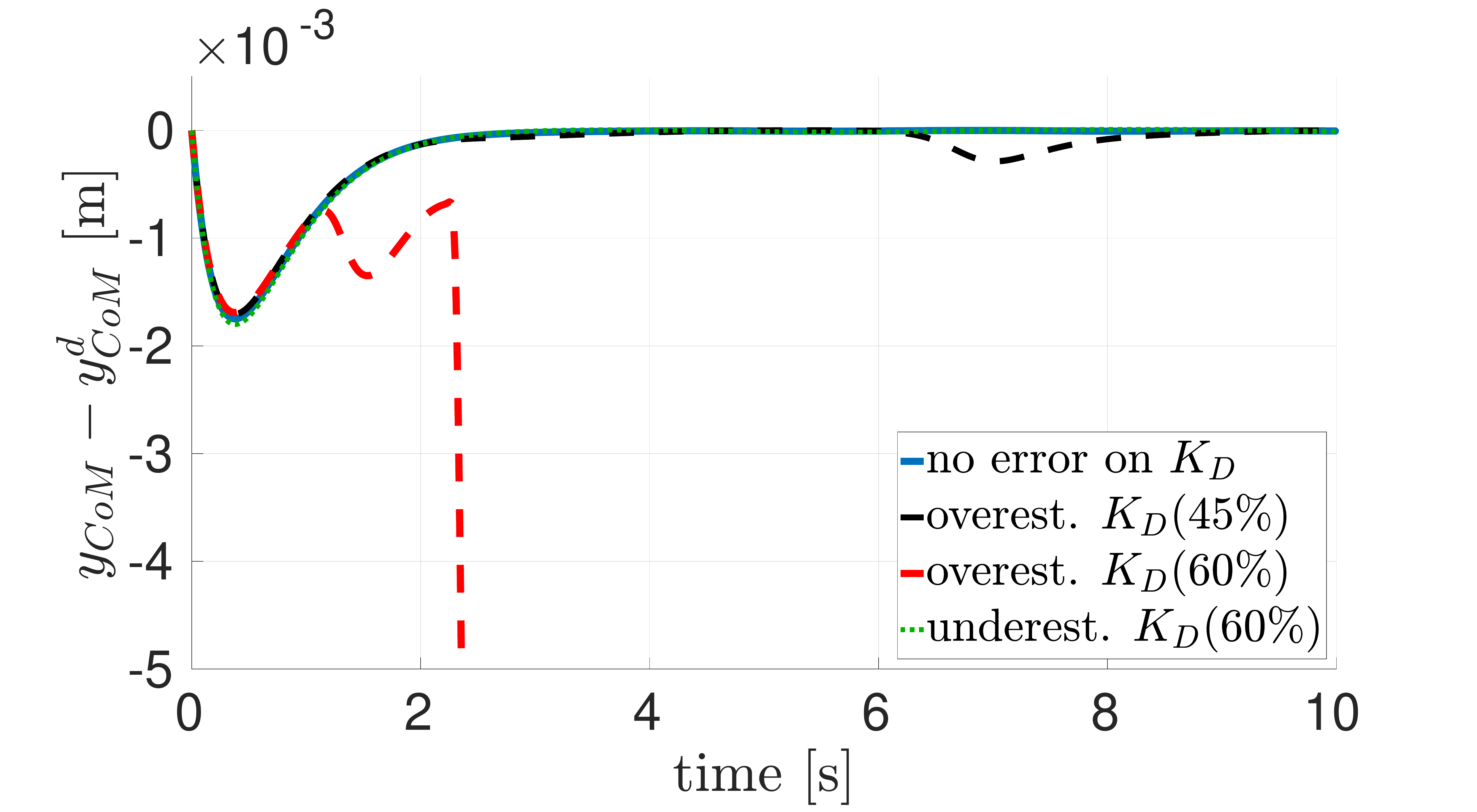}
   \caption{CoM error on lateral direction. If $K_D$ is overestimated by $60\%$, the closed-loop system is unstable.}
   \label{fig:sens}
\end{figure}

We evaluate the performances of control laws \eqref{inputTorquesSOTu}-\eqref{eq:tauMotors} for tracking a desired center of mass trajectory while the robot is balancing on one foot. The reference trajectory for the center of mass is a sinusoidal curve with amplitude of $1$ $[cm]$ and frequency of $0.2$ $[Hz]$ along the robot lateral direction.  
Furthermore, we perform an indicative sensitivity analysis on the tracking error in case of uncertainties on damping matrix $K_D$. This analysis holds for this particular task, but results may vary in case other movements are required. 
Figure~\ref{fig:sens} shows the tracking error on center of mass trajectory.
The thick blue line represents the center of mass error assuming perfect knowledge of system's damping. The dashed lines represent the error when the controller overestimates the real damping by $45\%$ and $60\%$, while the dotted line is obtained when $K_D$ is underestimated by $60\%$. Focus on the case when $K_D$ is overestimated by the controller, since it may lead to unstable behaviours. In particular, the center of mass error is still not increasing up to an error of $45\%$ on $K_D$, while an error of $60\%$ makes the system unstable.

\subsection{Comparisons between different control models}
\label{subsec:comparison}

Theoretically, one may apply the control law obtained from~\eqref{inputTorquesSOT}, which was developed under the assumption of stiff joints, when the system is, instead, powered by series elastic actuators, i.e. it is governed by~\eqref{eq:elasticSystem}. We here show that in this case, the closed loop system may have unstable behaviours because the control model neglects joint elasticity. In particular, the solution to the optimisation problem~\eqref{inputTorquesSOT} is a joint torque $\tau = \Gamma^{-1}\tau_m$ (see Remark \ref{remark1}).

\begin{figure}[t!]
   \centering
   \includegraphics[width=\columnwidth]{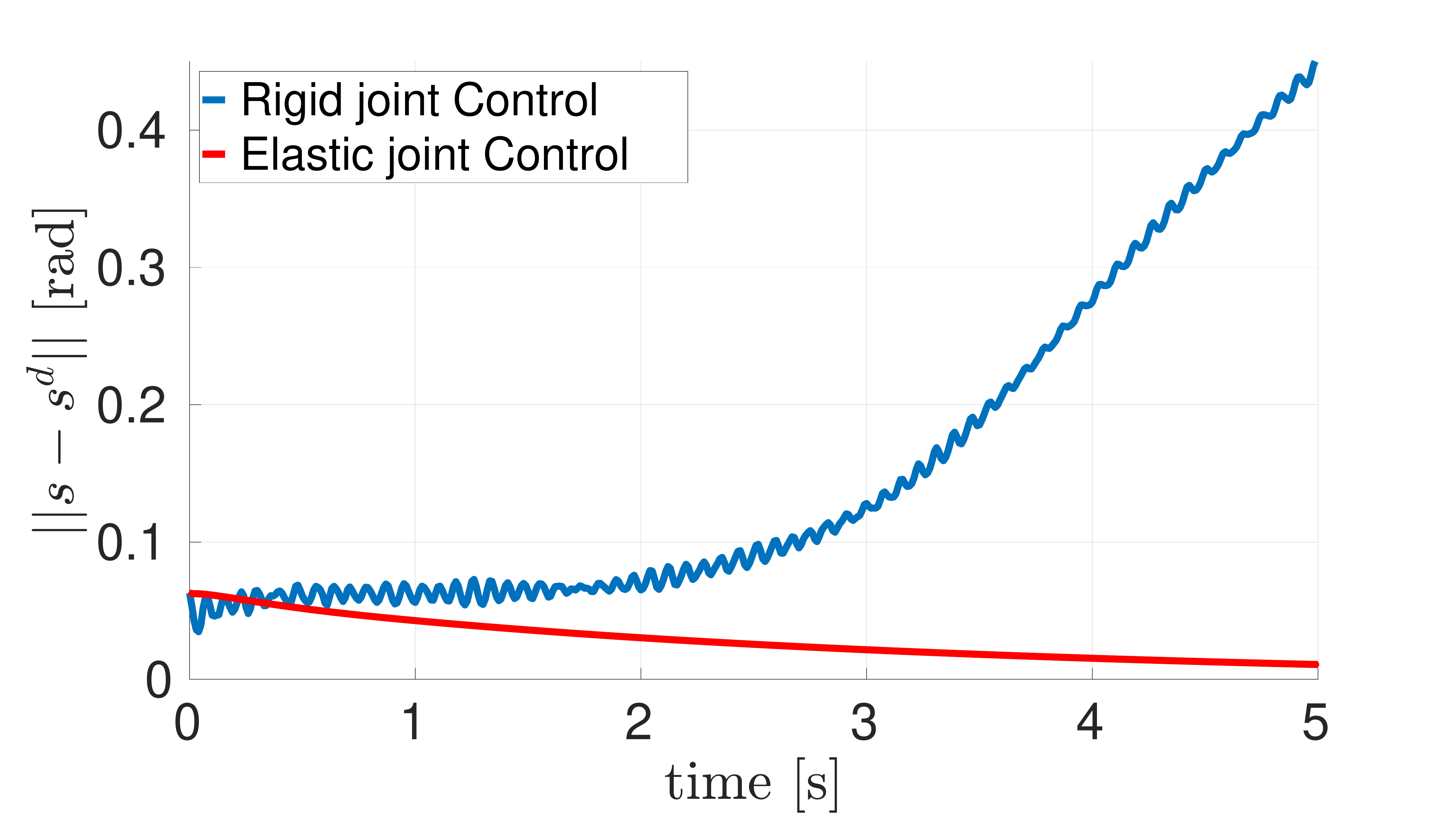}
   \caption{Comparison between control law \eqref{eq:forces}--\eqref{eq:torques}--\eqref{posturalNew} and \eqref{inputTorquesSOTu}-\eqref{eq:tauMotors} for controlling system \eqref{eq:elasticSystem}. The rigid joints control fails to stabilize the closed loop system about the reference position.}
   \label{fig:stability}
\end{figure}

To test the controller obtained from~\eqref{inputTorquesSOT} in the series-elastic actuator case, we impose a step response for the desired center of mass that translates into a   step of $1^\circ$ for all upper body joints (torso and arms). Figure~\ref{fig:stability} depicts  the norm of the joint position errors in both stiff and elastic control model. It is clear that the closed-loop system with the control \eqref{eq:forces}--\eqref{eq:torques}--\eqref{posturalNew} (blue line) is unstable, while the elastic joint controller (red line) ensures the convergence to the desired position.

\subsection{Effects of torque saturation}

\begin{figure}[t!]
   \centering
   \includegraphics[width=\columnwidth]{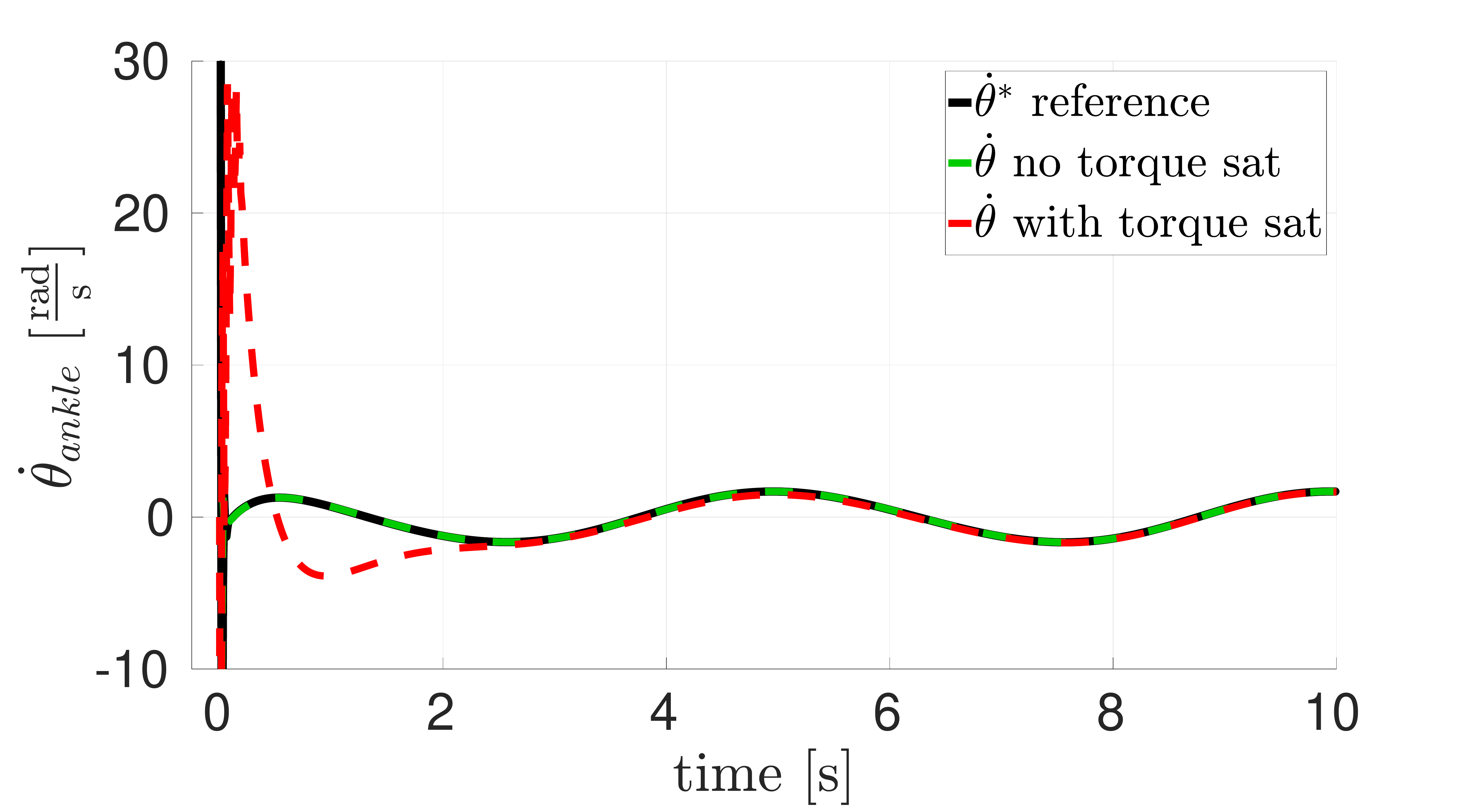}
   \caption{Convergence of motor velocity to $\dot{\theta}^*$ with and without torque saturation. Even in case of limited torques, the controller is able to stabilize the closed-loop system.}
   \label{fig:saturation}
\end{figure}

On the real robot, it may not be possible to achieve the desired motor velocities $\dot{\theta}_d$ because of limited motor torques, and the controller might fail to stabilize the closed-loop system.
To analyze the behaviour of elastic joint control in presence of limited motor torques, we add torque saturation to the simulation setup used in \ref{subsec:tracking}. We focused our attention on motor velocity of the stance foot ankle roll because this is the joint that requires the biggest torque for the given task. The maximum motor torque available is obtained through motors datasheet and it is $0.34$ $Nm$.
Figure \ref{fig:saturation} shows the effect of torques saturation on motor velocity: the black line is the reference velocity. To better visualize the results, we cut the initial peak (around $50$ $\frac{rad}{S}$). Dashed green line represents motor velocity without torque saturation, while the red line considers also torque saturation. In this second case, the convergence of $\dot{\theta}$ to $\dot{\theta}^*$ is slower, but stability is still retained.

\section{CONCLUSIONS}
\label{sec:conclusions}

This paper proposes a simple framework for extending the momentum based controllers developed for humanoid robots with stiff actuators to the case of series elastic actuators. It is based on the usage of motor velocities as a fictitious control input. Then,  fast convergence of the motor velocities to the desired values is obtained through feedback linearization of motors dynamics and (if necessary) high control gains for motor velocity error. Compared to other strategies, our control framework is robust against the feedforward terms usually needed by pure feedback linearisation techniques, and allows us to easily extend the momentum based controllers developed for rigid joints to the elastic joint case. 

In this paper, no experimental results are presented because series elastic actuators developed in \cite{Parmiggiani2012,Tisi2016} are about to be installed on real robot. Future work consists in validating the controller with series elastic actuators for both iCub version 2.5 \cite{Parmiggiani2012} and version 3 \cite{Tisi2016}.

\addtolength{\textheight}{0cm}     

\bibliographystyle{IEEEtran}
\bibliography{IEEEabrv,Biblio}

\end{document}